\documentclass[envcountsame,envcountsect,graybox]{svmult}

\usepackage{type1cm}        
\usepackage{makeidx}         
\usepackage{graphicx}        
\usepackage{multicol}        
\usepackage[bottom]{footmisc}

\usepackage{newtxtext}       %
\usepackage{newtxmath}       

\makeindex             

\usepackage{mdframed}
\usepackage{bm}
\usepackage{paralist}
\usepackage{tikz}
\usepackage{pgfplots}
\usepackage{pgfplotstable}

\usepackage[skip=2pt,font=scriptsize]{caption}
\usepackage{mathrsfs}

\newcommand{\cK}{\mathcal{K}}

\newcommand{\cR}{\mathcal{R}}

\newcommand{\cT}{\mathcal{T}}

\newcommand{\R}{\mathbb{R}}

\newcommand{\Dt}{{\Delta t}}

\DeclareMathOperator*{\Span}{span}

\makeatletter
\newcommand{\opnorm}{\@ifstar\@opnorms\@opnorm}
\newcommand{\@opnorms}[1]{%
  \left|\mkern-1.5mu\left|\mkern-1.5mu\left|
   #1
  \right|\mkern-1.5mu\right|\mkern-1.5mu\right|
}
\newcommand{\@opnorm}[2][]{%
  \mathopen{#1|\mkern-1.5mu#1|\mkern-1.5mu#1|}
  #2
  \mathclose{#1|\mkern-1.5mu#1|\mkern-1.5mu#1|}
}
\makeatother

\newcommand{\KU}[1]{{\textcolor{black}{#1}}}

\newcommand\OurAbstract{\KU{%
	Variational formulations of time-dependent PDEs in space and time yield $(d+1)$-dimensional problems to be solved numerically. This increases the number of unknowns as well as the storage amount. On the other hand, this approach enables adaptivity in space and time as well as model reduction w.r.t.\ both type of variables. In this paper, we show that matrix oriented techniques can significantly reduce the computational timings for solving the arising linear systems outperforming both time-stepping schemes and other solvers.
	}}

\begin{document}

\title*{Matrix oriented reduction of \KU{space-time Petrov-Galerkin variational} 
problems.}
\author{Julian Henning, Davide Palitta, Valeria Simoncini, and Karsten Urban}
\institute{Julian Henning, Karsten Urban  \at Ulm University, Institute for Numerical Mathematics,  
	Helmholtzstr.\ 20, 89081 Ulm (Germany), 
	\email{\{julian.henning,karsten.urban\}@uni-ulm.de}
	\and
	Davide Palitta \at
	Max Planck Institute for Dynamics of Complex Technical Systems, 
	Sandtorstr.\ 1, 39106 Magdeburg (Germany), 
	\email{palitta@mpi-magdeburg.mpg.de}
	\and
	Valeria Simoncini \at Universit\`{a} di Bologna, Centro AM$^2$, Dipartimento di Matematica, 
	Piazza di Porta S.\ Donato 5, 40127 Bologna (Italy), and IMATI-CNR, Pavia,
	\email{valeria.simoncini@unibo.it}
	}
\maketitle

\abstract*{\OurAbstract}

\abstract{\OurAbstract}

\section{Introduction}\label{sec:1}
Time-stepping schemes based upon variational semi-discretizations are the standard approach for the numerical solution of time-dependent partial differential equations (PDEs). Using a variational formulation in space and a subsequent discretization e.g.\ in terms of finite elements, one is left with an evolution problem in time. Standard finite difference techniques then yield a time-marching scheme, where a spatial problem needs to be solved in each time step.

Even though theoretical investigations on space-time variational formulations of PDEs 
{\color{black} have been} around for a long time, \cite{MR0247244}, it was seen prohibitive to treat the time as an additional variable for numerical simulations. In fact, if $\Omega\subset\R^d$ denotes the spatial domain, adding the time $t\in (0,T)=:I$ as an additional unknown results in a PDE on $\Omega_I:=I\times\Omega$ in Dimension $d+1$, which is costly both w.r.t.\ the amount of storage and the required computation time.

Also due to the increasing computing power the point of view has changed over the past years. 
In fact, being able to simulate problems for {\color{black} $d$ up to three until} about 20 years ago, 
adding another dimension seemed to be impossible. Nowadays, where we face high-dimensional problems (e.g.\ from quantum physics or finance) with $d\gg 100$, adding another dimension seems almost negligible. 
Another aspect to use space-time variational problems arose from model reduction of parameterized time-dependent PDEs. In fact treating both time and space as variables allows one to perform model reduction for space \emph{and} time, \cite{SpaceTimeUrbanPatera}. The time-stepping model reduction approach yields a time-marching scheme for a reduced spatial dimension but with the same number of time steps, \cite{MR3702337}.

In this paper, we address the question {\color{black} of} how to efficiently solve the linear systems arising from a (full) Petrov-Galerkin discretization of space-time variational formulations of time-dependent PDEs. 
{\color{black} It turns out that the involved coefficient matrices, combining space and time discretizations,
have a} tensorproduct structure, which allows {\color{black}us}
 to use {\color{black} more efficient matrix equations solvers than
 what can be done with the usual vector representation.} 

This paper is organized as follows: In \S\ref{sec:2} we review space-time variational formulations of some PDEs and describe corresponding Petrov-Galerkin discretizations as well as the arising linear systems in \S\ref{sec:3}; \S\ref{sec:4} is devoted to the description of the numerical schemes and \S\ref{sec:5} to numerical experiments, in particular the comparison with time-stepping schemes. 

\section{Space-Time Variational Formulation of PDEs}\label{sec:2}

\noindent\textbf{\KU{The Heat Equation.}} 
{\color{black}Let} $A:X\to X'$ be an elliptic operator on $X:=H^1_0(\Omega)$ associated to a coercive bilinear form $a:X\times X\to\R$, {\color{black} and $f\in L_2(I;X')$.}
We look for $u\in U:= H^1_{(0)}(I;X')\cap L_2(I;X)$ such that\footnote{$H^1_{(0)}(I;X'):= \{ w:I\to X':\, w\in H^1(I;X'), w(0)=0\}$, recall that $H^1(I;X')\hookrightarrow C(\bar{I};X')$.}
	$u_t + Au = f$, 
	$u(0)=0$, 
where 
homogeneous initial conditions are chosen only for convenience. The variational formulation then reads
\begin{align}\label{eq:var}
	\text{find } u\in U:\quad
	b(u,v) =  \langle f,v\rangle
	\quad\text{for all } v\in V,
\end{align}
where $V:= L_2(I;X)$, $b(u,v) := \int_0^T \int_\Omega u_t(t,x)\, v(t,x)\, dx\, dt + \int_0^T a(u(t), v(t))\, dt$ and $\langle f,v\rangle:= \int_0^T \int_\Omega f(t,x)\, v(t,x)\, dx\, dt$. The well-posedness is ensured by Ne\v{c}as' conditions, namely boundedness, injectivity and inf-sup condition of $b(\cdot,\cdot)$, \cite{MR1156075}.

\noindent\textbf{The Wave Equation.} Next, we consider an equation of wave type\KU{. Here, for $H:=L_2(\Omega)$, we view the operator $A$ as a mapping $A: \text{Dom}(A):=\{ \phi\in H: A\phi\in H\}\to H$, or $A: H\to \text{Dom}(A)'$. For $f\in L_2(I;H)$, we seek $u\in L_2(I;H)$ such that}
	$u_{t\!t} + Au = f$, 
	$u(0)=0$, $u_t(0)=0$, 
where we choose homogeneous initial conditions again only for convenience. 
In this case, it is not so obvious how to setup a well-posed variational form. It turns out that a very-weak setting is appropriate.
{\color{black} We formulate the problem as in (\ref{eq:var}) by using} $U:=L_2(I;H)$ as trial and $V:=\{ v\in L_2(I;H):\, v_{t\!t}+Av\in L_2(I;H),\, v(T)=v_t(T)=0\}$ as test space. Then, one can show that \eqref{eq:var} is well-posed for $b(u,v):=(u, v_{t\!t}+Av)_{L_2(I;H)}$ and $\langle f,v\rangle:= (f,v)_{L_2(I;H)}$ \KU{for fixed $T<\infty$}.

\section{Petrov-Galerkin Discretizations}\label{sec:3}
In order to determine a numerical approximation to the solution of a variational problem \eqref{eq:var}, one chooses finite-dimensional trial and test spaces, $U_\delta\subset U$, $V_\delta\subset V$, respectively. For convenience, we assume that their dimension is equal, i.e., $N_\delta:=\dim U_\delta=\dim V_\delta$. The Petrov-Galerkin method then reads
\begin{align}\label{eq:var-disc}
	\text{find } u_\delta\in U_\delta:\quad
	b(u_\delta,v_\delta) =  \langle f,v_\delta\rangle
	\quad\text{for all } v_\delta\in V_\delta.
\end{align}
As opposed to the coercive case, the well-posedness of \eqref{eq:var-disc} is not inherited from that of \eqref{eq:var}. In fact, the spaces $U_\delta$ and $V_\delta$ need to be appropriately chosen in the sense that the discrete inf-sup (or LBB --Ladyshenskaja-Babu\v{s}ka-Brezzi) condition holds, i.e., there exists an $\beta>0$ such that
$
	\inf_{u_\delta\in U_\delta} \sup_{v_\delta\in V_\delta} \frac{b(u_\delta,v_\delta)}{\| u_\delta\|_U\, \| v_\delta\|_V}
	\ge \beta >0,
$
where the crucial point is that $\beta\ne\beta_\delta$. The size of $\beta$ is also crucial for the error analysis, since it holds that $\| u-u_\delta\|_U \le \frac{1}{\beta} \inf_{w_\delta\in U_\delta}\| u-w_\delta\|_U$, \cite{MR1971217}.


\smallskip
\noindent\textbf{\KU{The Heat Equation}.} 
Starting with the temporal discretization, choose some integer $\KU{N_t}>1$ and set $\Dt:=T/\KU{N_t}$ resulting in a temporal triangulation 
$
	\cT_{\Dt}^\text{time}\equiv\{ t^{k-1}\equiv(k-1)\Delta t < t \le k\, \Dt \equiv t^k, 1\le k\le \KU{N_t}\}
	$ 
in time. Denote by $S_\Dt = \Span\{ \sigma^1,\ldots ,\sigma^{\KU{N_t}}\}$ piecewise linear finite elements on $I$, where $\sigma^k$ is the (interpolatory) hat-function with the nodes $t^{k-1}$, $t^k$ and $t^{k+1}$ (resp.\ truncated for $k\in\{0,\KU{N_t}\}$) and  $Q_\Dt=\Span\{ \tau^1,\ldots ,\tau^{\KU{N_t}}\}$ piecewise constant finite elements, where $\tau^k:=\chi_{I^k}$, the characteristic function on the temporal element $I^k:=(t^{k-1},t^k)$. 
For the spatial discretization, we choose any conformal $X_h=\Span\{ \phi_1,\ldots,\phi_{\KU{N_h}}\}\subset X$, e.g.\ piecewise linear finite elements. Then, we set $U_\delta := S_{\Delta t}\otimes X_{h}$, $V_\delta = Q_{\Delta t}\otimes X_{h}$, $\delta=(\Delta t,h)$. It can be shown that this yields LBB. Moreover, for $A=-\Delta$ and choosing the energy norm on $X$ as well as a slightly modified norm on $U$, one can even prove that $\beta=1$, \cite{SpaceTimeUrbanPatera}. Finally, we remark that this specific discretization  coincides with the Crank--Nicolson (CN) scheme if a trapezoidal approximation of the right-hand side temporal integration is used. Hence, we can later compare space-time Petrov-Galerkin numerical schemes with a CN time-stepping scheme. 

Finally, we detail the linear system of equations $B_\delta^T u_\delta = f_\delta$, where
\begin{align}
	[B_\delta]_{(k,i),(\ell,j)} 
		&= (\dot\sigma^k, \tau^\ell)_{L_2(I)}\, (\phi_i, \phi_j)_{L_2(\Omega)}
			+ (\sigma^k, \tau^\ell)_{L_2(I)}\, a(\phi_i, \phi_j),\label{eqn:par1} \\
	[f_\delta]_{(\ell,j)} 
		&= (f, \tau^\ell\otimes \phi_j)_{L_2(I;H)} \label{eqn:par2},
\end{align}
which means that we get a tensorproduct structure for the stiffness matrix $B_\delta= D_\Dt\otimes M_h + \KU{C_\Dt}\otimes A_h$, where the matrices are defined in an obvious manner. The right-hand side is not yet {\color{black} in a} 
tensorproduct structure. However, we can achieve that by determining an approximation
\begin{align}\label{eq:EIM}
	f(t,x) \approx \sum_{\KU{p}=1}^{\KU{P}} \vartheta_{\KU{p}}(t) f_{\KU{p}}(x) =: f^{\KU{P}}(t,x),
\end{align}
e.g.\ by the \emph{Empirical Interpolation Method} (EIM), \cite{MR2103208}. By choosing ${\KU{P}}$ sufficiently large, we can achieve any desired accuracy. Then, we get 
$[f_\delta^{\KU{P}}]_{(\ell,j)} 
		= \sum_{{\KU{p}}=1}^{\KU{P}}(\vartheta_{\KU{p}}, \tau^\ell)_{L_2(I)}\,  (f_{\KU{p}}, \phi_j)_{L_2(\Omega)}$, 
i.e., $f_\delta^{\KU{P}} = \sum_{{\KU{p}}=1}^{\KU{P}} h_{\KU{p}}\otimes g_{\KU{p}}$.

\smallskip
\noindent\textbf{The Wave Equation.} Constructing a stable pair of trial and test spaces for the wave equation is again a nontrivial task. Following an idea from \cite{JBSU18}, we first define the test space and construct the trial space in a second step in order to guarantee LBB\KU{, which, however, deteriorates with increasing $T$}.  Doing so, we set $R_\Dt := \Span\{ \varrho^1,\ldots ,\varrho^{\KU{N_t}}\}\subset H^2_{T}(I):=\{ \rho\in H^2(I): \rho(T)=\dot{\rho}(T)=0\}$, e.g.\ piecewise quadratic splines on $\cT_{\Dt}^\text{time}$. For space, we  choose any conformal $Z_h=\Span\{ \psi_1,\ldots,\psi_{\KU{N_h}}\}\subset H^2(\Omega)\cap H^1_0(\Omega)$, e.g.\ piecewise quadratic finite elements. Then, we define $V_\delta:= R_\Dt\otimes Z_h$, a tensor product space. 
The trial space $U_\delta$ is constructed by applying the adjoint PDE operator to each test basis function, i.e.
$
v_{k,i}:= \frac{d^2}{dt^2} \varrho^k(t) \psi_i(x) + A(\varrho^k(t) \psi_i(x))= \ddot{\varrho}^k(t) \psi_i(x) + \varrho^k(t) A\psi_i(x)
$.
We detail the arising linear system of equations starting with the stiffness matrix 
\begin{align*}
	[B_\delta]_{(k,i),(\ell,j)} 
		&= b(v_{k,i}, \varrho^\ell\kern-1pt\otimes\kern-1pt \psi_j)  
		= (\ddot{\varrho}^k\kern-1pt\otimes\kern-1pt \psi_i 
			+ \varrho^k \kern-1pt\otimes\kern-1pt A\psi_i,
				\ddot{\varrho}^\ell\kern-1pt\otimes\kern-1pt \psi_j 
			+ \varrho^\ell \kern-1pt\otimes\kern-1pt A\psi_j)_{L_2(I;H)}\\
		&= (\ddot{\varrho}^k, \ddot{\varrho}^\ell)_{L_2(I)}\,
			(\psi_i, \psi_j)_{L_2(\Omega)}
			+ (\ddot{\varrho}^k,\varrho^\ell)_{L_2(I)}\,
			 	 (\psi_i, A\psi_j)_{L_2(\Omega)}\\
		&\,\,+ (\varrho^k, \ddot{\varrho}^\ell)_{L_2(I)}\,
			(A\psi_i,\psi_j)_{L_2(\Omega)}
			+ (\varrho^k,\varrho^\ell)_{L_2(I)}\,
			(A\psi_i, A\psi_j)_{L_2(\Omega)},
\end{align*}
so that $B_\delta= Q_\Dt \otimes M_h + (D_\Dt+D_\Dt^T)\otimes A_h + M_\Dt\otimes Q_h$, again with obvious definitions of the matrices. For the right-hand side, we perform again an EIM-type approximation $f^{\KU{P}}(t,x)$. Then,
$
	[f^{\KU{P}}_\delta]_{(\ell,j)} 
		= \sum_{{\KU{p}}=1}^{\KU{P}} (\vartheta_{\KU{p}}\otimes f_{\KU{p}}, \varrho^\ell\otimes \psi_j)_{L_2(I;H)}
		= \sum_{{\KU{p}}=1}^{\KU{P}} (\vartheta_{\KU{p}}, \varrho^\ell)_{L_2(I)}\,
			 (f_{\KU{p}}, \psi_j)_{L_2(\Omega)}$,
so that the right-hand side has the same structure as in the first example. \KU{Due to the asymptotic behavior of the inf-sup-constant, we expect stability problems as $\Dt\to 0$, i.e., $N_t\to\infty$.}

\section{Efficient Numerical Methods for Tensorproduct Systems}\label{sec:4}


In both cases described above (and in fact also in space-time variational formulations of other PDEs), we obtain a (regular) linear system of the form 
\begin{equation}\label{eq:Kroneckersum}
	B_{\KU{\delta}}u_{\KU{\delta}}=f_{\KU{\delta}} 
	\quad\text{ with }\quad
	B_{\KU{\delta}}= \sum_{{\KU{p}}=1}^{{\KU{P}}_B} D_{\KU{p}}\otimes A_{\KU{p}}, \quad
	f_{\KU{\delta}} = \sum_{\ell=1}^{{\KU{P}}_f} h_\ell\otimes q_\ell, 
\end{equation}
\vskip -.05in
\noindent
where all involved matrices are sparse and (at least some of) the $A_q$ are s.p.d.
{\color{black} Recall that
$(D_{\KU{p}}\otimes A_{\KU{p}}) x = {\rm vec} (A_{\KU{p}} X D_{\KU{p}}^T)$, where {\rm vec} stacks the columns of
a given matrix one after the other, and $x={\rm vec}(X)$. We can thus rewrite the
system $B_\delta u_\delta = f_\delta$ in \eqref{eq:Kroneckersum} as the linear {\it matrix} equation 
$
	\sum_{{\KU{p}}=1}^{{\KU{P}}_B} A_{\KU{p}} U_\delta D_{\KU{p}}^T 
	= \sum_{\ell=1}^{{\KU{P}}_f} q_\ell h_\ell^T, 
$
with $u_\delta={\rm vec}(U_\delta)$. 
Matrix equations are receiving significant attention in the PDE context, due to the possibility
of maintaining the structural properties of the discretized problem, while limiting memory consumptions;
see \cite{Simoncini.survey.16}.
Under certain hypotheses, a large variety of discretization methodologies such as
finite differences, isogeometric analysis, spectral (element) methods, certain
finite element methods as well as various parametric numerical schemes
rely on tensor product spaces; see, e.g., 
{\cite{Bachmayretal.18,Dahmenetal.16,Griebel.Harbrecht:2011,KS11}.}
More recently, all-at-once time discretizations have shown an additional setting where tensor product 
approximations naturally arise; see, e.g., \KU{\cite{Palitta.19}} and references
therein.
Among the various computational strategies discussed in the literature \cite{Simoncini.survey.16},
here we focus on projection methods that reduce the original equation to a similar one, but of much smaller
dimension.

\smallskip
\noindent\textbf{Discretized \KU{Heat Equation}.}
The problem $B_\delta u_{\KU{\delta}} = f_\delta$ stemming from (\ref{eqn:par1},\ref{eqn:par2})
yields the following generalized Sylvester equation
\begin{align}\label{eq:Sylv_Heat}
	M_h U_\delta D_{\Delta t} + A_h U_\delta \KU{C}_{\Delta t}  = F\KU{_\delta}, \qquad {\rm with} 
	\quad 
	F\KU{_\delta} := [g_1, \ldots, g_{\KU{P}}] [h_1, \ldots, h_{\KU{P}}]^T.
\end{align}
The spatial stiffness and mass matrices $A_h$ and $M_h$ typically
have significantly larger dimensions \KU{$N_h$} than the time discretization matrices $D_{\Delta t}$, $N_{\Delta t}$, i.e., \KU{$N_t\ll N_h$}.
We therefore use a reduction method only for the space variables by projecting the problem onto
an appropriate space. A matrix Galerkin \KU{orthogonality} condition is then applied to obtain the solution: 
given $V_m \in {\mathbb R}^{\KU{N_h} \times k_m}$, $k_m\ll \KU{N_h}$, with orthonormal columns,
\KU{we} consider the approximation space range($V_m$) \KU{and} seek \KU{$Y_m\in\R^{k_m\times \KU{N_t}}$ such that } $U_{\delta,m}:= V_m Y_m \approx U_\delta$ \KU{and} the
residual $R_m := F\KU{_\delta} - (M_h U_{\delta,m} D_{\Delta t} + A_h U_{\delta,m} \KU{C}_{\Delta t})$ satisfies the Galerkin condition
$R_m \perp \text{range}(V_m)$. Imposing this orthogonality yields that $V_m^T R_m = 0$ is equivalent to 
$V_m^TF\KU{_\delta}V_m - (V_m^T M_h V_m) Y_\KU{m} D_{\Delta t} - (V_m^TA_h V_m) Y_\KU{m} \KU{C}_{\Delta t}=0$.
The resulting \KU{problem} is again a generalized Sylvester equation, but of much smaller size, therefore
Schur-decomposition oriented methods can cheaply be used, \cite[sec.4.2]{Simoncini.survey.16}, see
\cite{Simoncini.survey.16} for a discussion on projection methods \KU{as well as} their matrix and convergence properties.

For selecting $V_m$, let $F=F_1 F_2^T$ with $F_1$ \KU{having} full column rank.
Given the properties of $A_h$, $M_h$, we propose to employ the rational Krylov subspace 
$
\cR\cK_m := \textrm{range}([F_1, (A_h - \sigma_2 M_h)^{-1}M_h F_1, 
(A_h - \sigma_3 M_h)^{-1}M_h F_1, \ldots, (A_h - \sigma_m M_h)^{-1}M_h F_1]),
$
where the shifts $\sigma_s$ can be determined adaptively while the space is being
generated; see \cite{Simoncini.survey.16} for a description and references. The obtained
spaces are nested, $\cR\cK_m \subseteq \cR\cK_{m+1}$, therefore
the space can be expanded if the approximation is not sufficiently good. To include a
residual-based stopping criterion, the
residual norm can be computed in a cheap manner, see, e.g., \KU{\cite{Druskin.Simoncini.11,Palitta.19}}
for the technical details.

\noindent\textbf{Discretized Wave problem.} 
The problem $B_\delta u_{\KU{\delta}} = f_\delta$ now takes the matrix form
\begin{equation}\label{eq:3terms}
M_h U_\delta Q_{\Delta t}^T + A_h U_\delta (D_{\Delta t} + 
D_{\Delta t}^T) + Q_h U_\delta M_{\Delta t} = F_\delta.
\end{equation}
This three-term equation cannot be solved directly as before, therefore
we opt for using preconditioned GMRES on the vectorized equation. The preconditioner 
is given by the functional ${\cal P}: U \to M_h U Q_{\Delta t}^T + Q_h U M_{\Delta t}$,
corresponding to the discretized forth order operators, and exploits the matrix structure.
Hence, at the $k$th GMRES iteration 
we solve the generalized Sylvester equation $M_h W Q_{\Delta t}^T + Q_h W M_{\Delta t} = V_k$
where $V_k$ is such that $v_k={\rm vec}(V_k)$ is the previous basis vector. 
Since in this one-dimensional problem dimensions are limited, this matrix equation is solved
by explicitly diagonalizing the pairs 
$(Q_h, M_h)$ and $(Q_{\Delta t}, M_{\Delta t})$ \cite{Simoncini.survey.16}.

\section{Numerical Experiments}\label{sec:5}
{
In this section we show that the numerical solution of the linear system 
$B_\delta u_\delta=f_\delta$ can largely benefit from the exploitation of its 
Kronecker sum structure~\eqref{eq:Kroneckersum}. 
The performance of the all-at-once methods is compared in terms of both computational time 
and memory requirements. 
For the heat equation, We also document the performances of CN in terms of computational time. We are 
not aware of any variant of CN that is able to exploit the low-rank structure of the underlying 
problem and we thus employ the classical CN scheme. Such implementation leads to 
large running times and a storage demand that is always equal to 
\KU{$N_t\cdot N_h$} as the full $U_\delta$ is allocated.

The tolerance of the final relative residual norm is set to $10^{-8}$ and in the following tables
we also report the number of iterations needed to achieve such accuracy and the numerical
rank of the computed solution. 
All results were obtained with Matlab R2017b on a machine with 2.1 GHz processors and 192GB
of RAM.
}

\begin{example}\label{Ex.1}
\emph{(The heat equation)} We consider the equation  on 
the cube $\Omega = (-1,1)^3$  with homogeneous Dirichlet boundary conditions and the time 
interval $I \coloneqq (0, 10)$ with initial conditions  $u(0,x,y,z) \equiv 0$.  
The right-hand side is $f(t,x,y,z) := 10 \sin(t)t \cos(\frac{\pi}{2}x )\cos(\frac{\pi}{2}y )\cos(\frac{\pi}{2}z )$
and its discretized version is thus low rank.
For discretization in space, linear finite elements were chosen, leading
to the discretized generalized Sylvester equation in (\ref{eq:Sylv_Heat}). 
We compare the performance of the Galerkin projection method based upon rational Krylov
spaces described in \S\ref{sec:4} (denoted RKSM) with that of a low-rank version of
preconditioned GMRES (denoted LR-FGMRES-RKSM). For the latter, 
the preconditioner is chosen as a fixed (five) number of iterations 
of the rational Krylov Galerkin projection method. 
Since the preconditioner is a non-linear operator, a flexible variant of GMRES is used.
We refer the reader to \cite{Palitta.19} for more details.
The results are displayed in Table \ref{tab1.1}. 

\begin{table}[htb]
 \centering
 \caption{Results for Example \ref{Ex.1}: different values of $N_h$ and $N_t$. Memory allocations 
for RKSM and LR-FGMRES+RKSM are given by $\mu_{\text{mem}}\cdot(N_h+N_t)$.\label{tab1.1}}
\begin{tabular}{| r r|rrrr |rrrr|r| }
\hline
  & & \multicolumn{4}{c|}{RKSM} & \multicolumn{4}{c|}{LR-FGMRES+RKSM}& \multicolumn{1}{c|}{CN}\\
   \KU{$N_h$} & $N_t$ & Its & $\mu_{\text{mem}}$ & $\text{rank}(U_\delta)$ & Time (s) & Its & $\mu_{\text{mem}}$ & $\text{rank}(U_\delta)$ & Time (s)&   Time (s) \\
 \hline 
 41300 & 100 & 16& 17 & 10& 88.74
 & 4 & 83& 13 & 143.13  & 296.16  \\
       
       & 300 & 13& 14 & 9& 68.86 & 4 & 74& 10 & 130.27 & 871.38  \\
       
       
       & 500 & 13& 14 & 9& 65.88 & 4 & 75& 11 &134.73 & 1468.40  \\
\hline
347361 & 100 & 16& 17 & 9& 3144.02
 & 4 & 71& 10 & 4793.63  & 13805.09  \\
       
       & 300 & 14& 15 & 9& 2673.83 & 4 & 78& 9 & 4780.46 & 41701.10  \\
       
       
       & 500 & 14& 15 & 9& 2699.98 & 4 & 80& 9 &4827.13 & 70044.52  \\
\hline
\end{tabular}
 \end{table}

 The CN method leads to an excessive workload compared with the all-at-once approaches
for all considered values of $N_h$ and $N_t$, with the computational time growing linearly 
with the number of time steps $N_t$. The performance of the other methods is independent
of the time discretization, and it only depends on the spatial component of the overall discrete operator.
In fact, spatial mesh independence seems to also be achieved.
%
 
 \end{example}
 
\begin{example}\label{Ex.2}
\emph{(The wave equation)} 
\KU{We consider the wave problem with $A = -\Delta$ on $\Omega = (0,1)$ with homogeneous Dirichlet boundary conditions and $I \coloneqq (0, 1)$.  Setting $f(t,x) := \sin(2 \pi x) +  4 \pi^2 t^2 \sin(2 \pi x)$  yields the analytical solution $u(t,x) = t^2 \sin(2 \pi x)$. We choose cubic B-Splines for the discretization in space and time. }
The discretized problem thus leads to the matrix equation in (\ref{eq:3terms}).
In Figure~\ref{tab2.1} we report some our preliminary results. Note, that our above discretization does not yield an equivalent time-stepping scheme with which we could do comparisons.

\begin{figure}   
\caption{Example~\ref{Ex.2}. Left: Results for different values of $N_h$ and $N_t$. Right: Relative residual norm history for some values of $N_h$ and $N_t$. \vspace*{-3mm}} 
  \label{tab2.1}
	\begin{center}
	\begin{tabular}{cc}
		{	
 		\setlength\tabcolsep{0.3pt}
 		\begin{tabular}{ r r rrrr }
 		\vspace{-4.7 cm}
 		&&&&&\\
\hline
  \multicolumn{2}{|c|}{} & \multicolumn{3}{c}{GMRES+LYAP} & \multicolumn{1}{|c|}{{\tt backslash}}\\
    \multicolumn{1}{|r}{\KU{$N_h$}} & \multicolumn{1}{r|}{\KU{$N_t$}} & Its. & $\text{rank}(U_\delta)$ & Time (s) &\multicolumn{1}{|r|}{Time (s)} \\
 \hline 
 \multicolumn{1}{|r}{256} & \multicolumn{1}{r|}{256}  & 16 & 13 & 0.21 & \multicolumn{1}{|r|}{1.17} \\

 \multicolumn{1}{|r}{}  & \multicolumn{1}{r|}{512} & 36 & 35 & 1.51 & \multicolumn{1}{|r|}{2.39} \\

 \multicolumn{1}{|r}{}    & \multicolumn{1}{r|}{1024} & 81 & 74 & 20.97 & \multicolumn{1}{|r|}{12.61} \\
    \hline
    
    \multicolumn{1}{|r}{512} & \multicolumn{1}{r|}{256} & 26 & 31 & 0.61 & \multicolumn{1}{|r|}{2.30} \\
    
 \multicolumn{1}{|r}{}    & \multicolumn{1}{r|}{512} & 40 & 43 &  2.64 & \multicolumn{1}{|r|}{5.09} \\
  
  \multicolumn{1}{|r}{}   & \multicolumn{1}{r|}{1024} & 81 & 74 & 20.97 & \multicolumn{1}{|r|}{12.61}\\
 \hline
 
 \multicolumn{1}{|r}{1024} & \multicolumn{1}{r|}{256} & 50 & 59 & 3.55 & \multicolumn{1}{|r|}{4.82} \\
 
 \multicolumn{1}{|r}{}   & \multicolumn{1}{r|}{512} & 68 & 72 & 10.12 & \multicolumn{1}{|r|}{11.13} \\

 \multicolumn{1}{|r}{}   & \multicolumn{1}{r|}{1024} & 102 & 92 & 54.15 & \multicolumn{1}{|r|}{24.28} \\    
  \hline
\end{tabular}
}		&
 		{  
 		\includegraphics[scale=0.401]{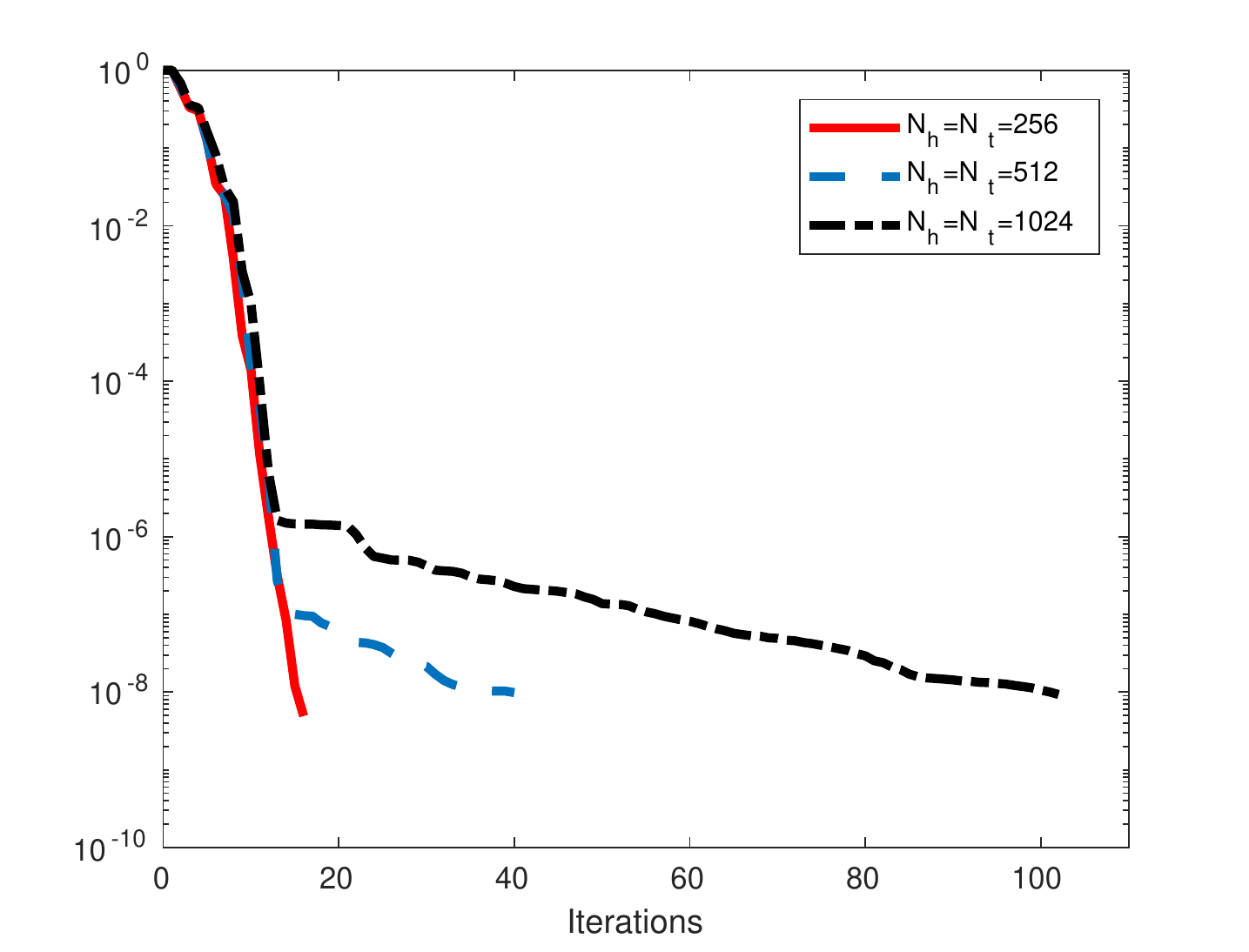}
 		
}\\
	\end{tabular}
	\end{center}
\end{figure}
\vspace{-1cm}

The table on the left shows that the performances of our preconditioned scheme are quite good for small values of $N_h$ and $N_t$. Indeed, in this case, the preconditioner manages to drastically reduce the number of iterations needed to converge so that GMRES+LYAP 
turns out to be faster than the Matlab solver {\tt backslash} applied to the solution of the linear system $B_\delta u_\delta=f_\delta$, in spite of the 1D nature (in space) of the problem. 
However, the effectiveness of the adopted preconditioner worsens by increasing the number of degrees of freedom. This is due to a dramatic increment in the condition number of the coefficient matrices (see the discussion at the end of \S\ref{sec:3}) that causes an abrupt very slow decrement (almost stagnation) in 
the GMRES residual at the level that seems to be related to the conditioning of 
the involved matrices, see Fig.\ \ref{tab2.1} (right). 
As it is, the problem associated with handling this ill-conditioning in the algebraic
equation is crucial for the overall solver performance, and will be the topic of future works.
Alternatively, one may try to directly address the solution of the multiterm 
matrix equation (\ref{eq:3terms})  as it is 
done in \cite{Powell2017} for certain stochastic PDEs. 
Projection-based solvers may be very well-suited for our framework as the solution $U_\delta$ turns out to be low-rank.
However, such an approach needs the design of ad hoc approximation spaces and we plan to pursue
this challenging research direction in the near future.

\end{example}

%
%
%
%

\bibliographystyle{abbrv}
\bibliography{ST_references}

\end{document}